\documentclass[amstex,12pt, amssymb]{article}

\usepackage{mathtext}
\usepackage[cp1251]{inputenc}
\usepackage[T2A]{fontenc}
\usepackage[dvips]{graphicx}
\usepackage{amsmath}
\usepackage{amssymb}
\usepackage{amsxtra}
\usepackage{latexsym}
\usepackage{ifthen}

\newcounter{lemma}[section]

\newcounter{corollary}[section]

\newcounter{remark}[section]

\newcounter{theorem}[section]

\newcounter{proposition}[section]

\numberwithin{equation}{section}

\begin{document}

\markboth{\centerline{E.~SEVOST'YANOV, S.~SKVORTSOV}}
{\centerline{ON RING HOMEOMORPHISMS... }}

\def\cc{\setcounter{equation}{0}
\setcounter{figure}{0}\setcounter{table}{0}}

\overfullrule=0pt

%\normalsize\large

\author{{E.~SEVOST'YANOV, S.~SKVORTSOV}\\}

\title{
{\bf ON RING HOMEOMORPHISMS WITH INVERSE MODULUS CONDITIONS}}

\date{\today}
\maketitle

%\large
\begin{abstract} We consider a class of so-called ring $Q$-mappings that are a generalization of quasiconformal
mappings. Theorems on the local behavior of inverse maps of this
class are obtained. Under certain conditions, we also investigated
the behavior of families of these mappings in the closure of the
domain.
\end{abstract}


\section{Introduction} The paper is devoted to the study of quasiconformal mappings and mappings with finite distortion,
actively studied recently (see, e.g., \cite{MRV$_2$}--\cite{Va},
cf.~\cite{Sev$_3$}.

\medskip
Let $M$ be the modulus of family of paths (see \cite{Va}), and
$dm(x)$ corresponds to the Lebesgue measure in ${\Bbb R}^n,$
$n\geqslant 2.$ Given sets $E, F$ and $D$ in $\overline{{\Bbb
R}^n}={\Bbb R}^n\cup \{\infty\},$ $\Gamma(E, F, D)$ denotes the
family of all paths $\gamma:[0, 1]\rightarrow \overline{{\Bbb R}^n}$
such that $\gamma(0)\in E,$ $\gamma(1)\in F$ and $\gamma(t)\in D$
for all $t\in (0, 1).$ In what follows, the boundary and the closure
of the set are understood in the sense of $\overline{{\Bbb R}^n}.$
Let $x_0\in\overline{D},$ $x_0\ne\infty,$
$$S(x_0,r) = \{
x\,\in\,{\Bbb R}^n : |x-x_0|=r\}\,, S_i=S(x_0, r_i)\,,\quad
i=1,2\,,$$
$$A=A(x_0, r_1, r_2)=\{ x\,\in\,{\Bbb R}^n : r_1<|x-x_0|<r_2\}\,.$$
Let $Q:{\Bbb R}^n\rightarrow {\Bbb R}^n$ be a Lebesgue measurable
function, $Q(x)\equiv 0$ in ${\Bbb R}^n\setminus D.$ A mapping
$f:D\rightarrow \overline{{\Bbb R}^n}$ is said to be {\it ring
$Q$-mapping at the point $x_0$}, if
\begin{equation} \label{eq2*!}
M(f(\Gamma(S_1, S_2, D)))\leqslant \int\limits_{A\cap D} Q(x)\cdot
\eta^n (|x-x_0|)\, dm(x)
\end{equation}
holds for each $0<r_1<r_2<d_0:=\sup\limits_{x\in D}|x-x_0|,$ and any
measurable function
$\eta:(r_1, r_2)\rightarrow [0, \infty]$ such that
\begin{equation}\label{eq8B}
\int\limits_{r_1}^{r_2}\eta(r)dr\geqslant 1
\end{equation}
(see, e.g.,~\cite{RSY$_1$}, cf.~\cite{MRSY$_4$}).
We say that $f$ is ring $Q$-mapping at $\infty,$ if $f(x/|x|^2)$ is
a ring $Q(x/|x|^2)$-mapping at the origin. A mapping $f:D\rightarrow
\overline{{\Bbb R}^n}$ is said to be {\it ring $Q$-mapping in $E$},
$E\subseteq \overline{D},$ if (\ref{eq2*!}) holds for every $x_0\in
E.$ If, in addition, $f$ is a homeomorphism, we say that $f$ is {\it
ring $Q$-homeomorphism in $E$}.

\medskip
The main definitions and notations used below can be found in
monographs \cite{Va} and \cite{MRSY}, and therefore are omitted.
Recall that the domain $D\subset {\Bbb R}^n$ is called {\it locally
connected at the point} $x_0\in\partial D,$ if for every
neighborhood $U$ of a point $x_0$ there is a neighborhood $V\subset
U$ of a point $x_0$ such that $V\cap D$ is connected. The domain $D$
is locally connected in $\partial D,$ if $D$ is locally connected at
every point $x_0\in\partial D.$ The boundary of $D$ is called {\it
weakly flat} at a point $x_0\in
\partial D,$ if for every $P>0$ and every neighborhood $U$
of the point $x_0,$ there is a neighborhood $V\subset U$ of $x_0$
such that $M(\Gamma(E, F, D))>P$ for all continua $E, F\subset D,$
intersecting $\partial U$ and $\partial V.$ The boundary of the
domain $D$ is weakly flat, if it is weakly flat at every point of
$\partial D.$

\medskip
Let $(X, d)$ and $\left(X^{{\,\prime}}, d^{\,\prime}\right)$ be
metric spaces with distances $d$ and $d^{\,\prime},$ respectively. A
family $\frak{G}$ of mappings $g:X^{\,\prime}\rightarrow X$  is said
to be {\it equicontinuous at a point} $y_0\in X^{\,\prime},$ if for
every $\varepsilon>0$ there is $\delta=\delta(\varepsilon, y_0)>0$
such that $d(g(y), g(y_0))<\varepsilon$ for all $g\in \frak{G}$ and
$y\in X^{\,\prime}$ with $d^{\,\prime}(y , y_0)<\delta$. The family
$\frak{G}$ is {\it equicontinuous} if $\frak{G}$ is equicontinuous
at every point $y_0\in X^{\,\prime}.$ In what follows, we consider
that $X=D,$ where $D$ is a bounded domain in ${\Bbb R}^n,$ and $d(x,
y)=|x-y|.$ Besides that, $X^{\,\prime}=D^{\,\prime}$ or
$X^{\,\prime}=\overline{D^{\,\prime}}$ depending on the context,
where $D^{\,\prime}$ is a domain in $\overline{{\Bbb R}^n},$ and
$d^{\,\prime}(x, y)=h(x, y),$
$$h(x,y)=\frac{|x-y|}{\sqrt{1+{|x|}^2} \sqrt{1+{|y|}^2}}\,,\quad x\ne
\infty\ne y\,, \quad\,h(x,\infty)=\frac{1}{\sqrt{1+{|x|}^2}}\,.
$$
Given a set $E\subset\overline{{\Bbb R}^n},$ we put
\begin{equation}\label{eq9C}
h(E):=\sup\limits_{x,y\in E}h(x, y)\,,
\end{equation}
where $h(E)$ is called the {\it chordal (spherical)} diameter of
$E.$
Given domains $D, D^{\,\prime}\subset \overline{{\Bbb R}^n},$
$n\geqslant 2,$ and Lebesgue measurable function $Q:{\Bbb
R}^n\rightarrow [0, \infty],$ $Q(x)\equiv 0$ for $x\not\in D,$
denote ${\frak R}_Q(D, D^{\,\prime})$ the family of all
homeomorphisms $g:D^{\,\prime}\rightarrow D$ of $D^{\,\prime}$ onto
$D$ such that $f=g^{\,-1}$ is ring $Q$-homeomorphism in $D.$ The
following assertion is valid.

\medskip
\begin{theorem}\label{th1}
{\sl Let $D$ be a bounded domain in ${\Bbb R}^n.$ If $Q\in L^1(D),$
then the family ${\frak R}_Q(D, D^{\,\prime})$ is equicontinuous in
$D^{\,\prime}.$ }
\end{theorem}

\medskip
Given $\delta>0,$ domains $D$ and $D^{\,\prime}\subset
\overline{{\Bbb R}^n},$ $n\geqslant 2,$ a continuum $A\subset D$ and
Lebesgue measurable function $Q:{\Bbb R}^n\rightarrow [0, \infty],$
$Q(x)\equiv 0$ for $x\not\in D,$ denote ${\frak S}_{\delta, A, Q
}(D, D^{\,\prime})$ the family of all homeomorphisms
$g:D^{\,\prime}\rightarrow D$ of $D^{\,\prime}$ onto $D$ such that
$f=g^{\,-1}$ is ring $Q$-homeomorphism in $\overline{D},$ wherein
$$h(f(A)):=\sup\limits_{x, y\in f(A)}h(x, y)\geqslant\delta\,.$$ The following
assertion is valid.

\medskip
\begin{theorem}\label{th2}
{\sl Let $D$ be a bounded domain in ${\Bbb R}^n.$ Suppose that $D$
is locally connected on the boundary, $\partial D^{\,\prime}$ is
weakly flat, and any connected component of $\partial D^{\,\prime}$
does not degenerate to a point. If $Q\in L^1(D),$ then each mapping
$g\in {\frak S}_{\delta, A, Q }(D, D^{\,\prime})$ have a continuous
extension $\overline{g}:\overline{D^{\,\prime}}\rightarrow
\overline{D},$ $\overline{g}|_{D^{\,\prime}}=g$ and
$\overline{g}(\overline{D^{\,\prime}})=\overline{D}.$ Moreover, the
family ${\frak S}_{\delta, A, Q }(\overline{D},
\overline{D^{\,\prime}}),$ consisting of all extended mappings
$\overline{g}:\overline{D^{\,\prime}}\rightarrow \overline{D},$ is
equicontinuous in $\overline{D^{\,\prime}}.$ }
\end{theorem}

\section{Preliminaries}
First of all, we establish two elementary statements that play an
important role in the proof of the main results. Let $I$ be an open,
closed or half-open interval in ${\Bbb R}.$ As usual, for a path
$\gamma: I\rightarrow {\Bbb R}^n,$ we set
$$|\gamma|=\{x\in {\Bbb R}^n: \exists\, t\in [a, b]:
\gamma(t)=x\}\,,$$
wherein, $|\gamma|$ is called {\it locus (image) of the path}
$\gamma.$ We say that the path $\gamma$ lies in the domain $D,$ if
$|\gamma|\subset D.$ Besides that, we say that paths $\gamma_1$ and
$\gamma_2$ are disjoint, if their loci do not intersect. The path
$\gamma:I\rightarrow {\Bbb R}^n$ is called {\it Jordan arc}, if
$\gamma$ is a homeomorphism of $I.$ The following (almost obvious)
assertion is valid.

\medskip
\begin{lemma}\label{lem1}{\sl\,
Let $D$ be a domain in ${\Bbb R}^n,$ $n\geqslant 2,$ locally
connected on its boundary.  Then any two pairs of different points
$a\in D, b\in \overline{D},$ и $c\in D, d\in \overline{D}$ can be
joined by disjoint paths $\gamma_1:[0, 1]\rightarrow \overline{D}$
and $\gamma_2:[0, 1]\rightarrow \overline{D},$  so, that
$\gamma_i(t)\in D$ for all $t\in (0, 1),$ $i=1,2,$ $\gamma_1(0)=a,$
$\gamma_1(1)=b,$ $\gamma_2(0)=c,$ $\gamma_2(1)=d.$}
\end{lemma}

\medskip
The following lemma shows that inner points of each domain are
<<weakly flat>>.

\medskip
\begin{lemma}\label{lem2}
{\sl\, Let $D$ be a domain in $\overline{{\Bbb R}^n},$ $n\geqslant
2,$ and $x_0\in D.$ Then, for every $P>0$ and for any neighborhood
$U$ of the point $x_0$ there is a neighborhood $V\subset U$ of the
same point, such that $M(\Gamma(E, F, D))>P$ for any continua $E,
F\subset D,$ intersecting $\partial U$ and $\partial V.$}
\end{lemma}

\section{Proof of Theorem~\ref{th1}}
\begin{figure}[h]
\centerline{\includegraphics[scale=0.6]{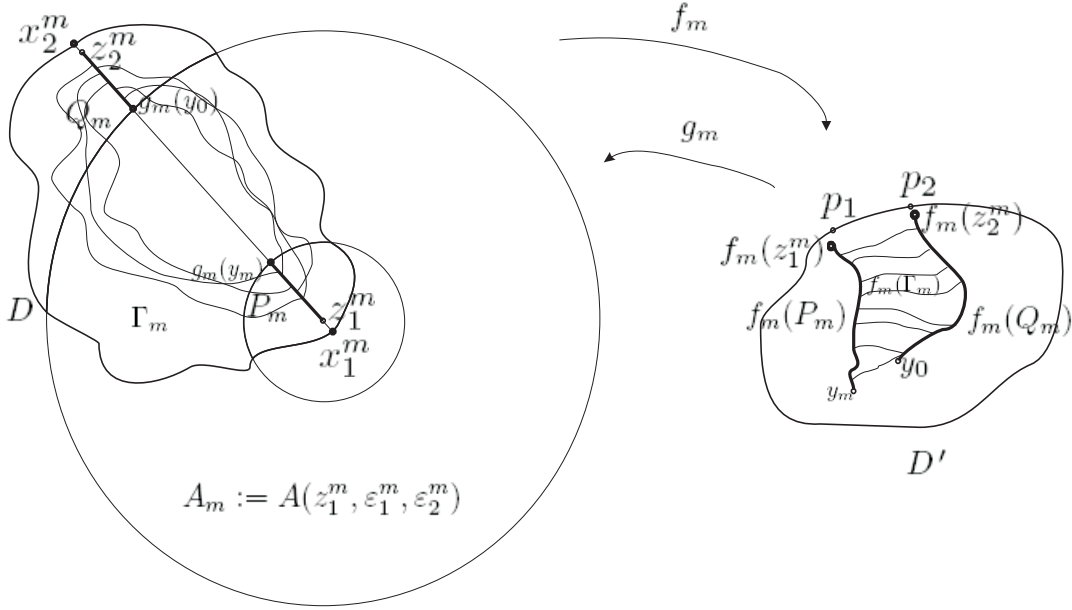}} \caption{To
the proof of the theorem~\ref{th1}}\label{fig2}
\end{figure}

We prove the theorem~\ref{th1} by contradiction. Suppose, the family
${\frak R}_Q(D, D^{\,\prime})$ is not equ\-i\-con\-ti\-nuous at some
point $y_0\in D^{\,\prime},$ in other words, there are $y_0\in
D^{\,\prime}$ and $\varepsilon_0>0,$ such that for any $m\in {\Bbb
N}$ there exists $y_m\in D^{\,\prime},$ $h(y_m, y_0)<1/m,$ and a
homeomorphism $g_m\in{\frak R}_Q(D, D^{\,\prime}),$ for which
\begin{equation}\label{eq13***}
|g_m(y_m)-g_m(y_0)|\geqslant \varepsilon_0\,.
\end{equation}
Let us consider the straight line
$$r=r_m(t)=g_m(y_0)+(g_m(y_m)-g_m(y_0))t,\quad-\infty<t<\infty,$$
passing through points $g_m(y_m)$ and $g_m(y_0)$
(see~Figure~\ref{fig2}).
Since $D$ is bonded, by \cite[Theorem~1.I, Ch.~5, \S\, 46]{Ku},
$r=r_m(t)$ intersects $\partial D$ for some $t\geqslant 1.$ Thus,
there exists $t_1^m\geqslant 1$ such that $r_m(t^m_1)=x^m_1\in
\partial D.$ Without loss of generality, we may assume that $r_m(t)\in
D$ for all $t\in [1, t^m_1).$ Now, the segment
$\gamma^m_1(t)=g_m(y_0)+(g_m(y_m)-g_m(y_0))t,$ $t\in [1, t^m_1],$
belongs to $D$ for all $t\in [1, t^m_1),$
$\gamma^m_1(t^m_1)=x^m_1\in
\partial D$ and $\gamma^m_1(1)=g_m(y_m).$
Similarly, there are $t^m_2<0$ and a segment
$\gamma^m_2(t)=g_m(y_0)+(g_m(y_m)-g_m(y_0))t,$ $t\in [t^m_2, 0],$
such that $\gamma^m_2(t^m_2)=x^m_2\in
\partial D,$ $\gamma^m_2(0)=g_m(y_0)$ and $\gamma^m_2(t)\in D$
for all $t\in (t^m_2, 0].$ Set $f_m:=g_m^{\,-1}$ and fix $m\in {\Bbb
N}.$ Since $f_m$ is a homeomorphism, $C(f_m, x_1^m)$ and $C(f_m,
x_2^m)$ belong to $\partial D^{\,\prime},$ where, as usually,
$$C(f, x):=\{y\in \overline{{\Bbb R}^n}:\exists\,x_k\in D: x_k\stackrel{k\rightarrow\infty}{\rightarrow}x, f(x_k)
\stackrel{k\rightarrow\infty}{\rightarrow}y\}\,,$$
see, e.g.,~\cite[Proposition~13.5]{MRSY}). Consequently, there is a
point $z^m_1\in D\cap |\gamma^m_1|$ such that ${\rm
dist}\,(f_m(z_1^m),
\partial D^{\,\prime})<1/m.$ Since $\overline{{\Bbb R}^n}$ is compact,
we can consider that $f_m(z_1^m)\rightarrow p_1\in \partial
D^{\,\prime}$ as $m\rightarrow\infty.$ Similarly, there is a
sequence $z^m_2\in D\cap |\gamma^m_2|$ such that ${\rm
dist}\,(f_m(z^m_2),
\partial D^{\,\prime})<1/m$ and
$f_m(z^m_2)\rightarrow p_2\in
\partial D^{\,\prime}$ as $m\rightarrow\infty.$

Let $P_m$ be the part of the interval $\gamma_1^m,$ enclosed between
the points $g_m(y_m)$ and $z^m_1,$ and $Q_m$ be the part of the
interval $\gamma_2^m,$ enclosed between the points $g_m(y_0)$ and
$z^m_2.$
Consider
$$A_m:=A(z_1^m, \varepsilon_1^m, \varepsilon_2^m)=\{x\in {\Bbb R}^n: \varepsilon_1^m<|x-z_1^m|<\varepsilon_2^m\}\,,$$
where
$$\varepsilon_1^m:=|g_m(y_m)- z_1^m|,\quad
\varepsilon_2^m:=|g_m(y_0)-z_1^m|\,.$$
Let $\Gamma_m=\Gamma(P_m, Q_m, D).$ Let us to prove that
\begin{equation}\label{eq1A}
\Gamma_m>\Gamma(S(z_1^m, \varepsilon_1^m), S(z_1^m,
\varepsilon_2^m), A_m)\,.
\end{equation}

Indeed, let $\gamma\in \Gamma_m,$ i.e., $\gamma=\gamma(s):[0,
1]\rightarrow {\Bbb R}^n,$ $\gamma(0)\in P_m,$ $\gamma(1)\in Q_m$
and $\gamma(s)\in D$ for $0<s<1.$ Let $q_m>1$ be a number, such that
$$z_1^m=g_m(y_0)+(g_m(y_m)-g_m(y_0))q_m\,.$$ Since $\gamma(0)\in P_m,$
there exists $1<t_m<q_m$ such that
$$\gamma(0)=g_m(y_0)+(g_m(y_m)-g_m(y_0))t_m\,.$$ Thus,
$$|\gamma(0)-z_1^m|=|(g_m(y_m)-g_m(y_0))(q_m-t_m)|< $$
\begin{equation}\label{eq2A}<
|(g_m(y_m)-g_m(y_0))(q_m-1)|=|(g_m(y_m)-g_m(y_0))q_m+g_m(y_0)-g_m(y_m))|=\end{equation}
$$=|g_m(y_m)-z_1^m|=\varepsilon_1^m\,.$$

From other hand, since $\gamma(1)\in Q_m,$ there exists $p_m<0$ such
that $$\gamma(1)=g_m(y_0)+(g_m(y_m)-g_m(y_0))p_m\,.$$ Now
$$
|\gamma(1)-z_1^m|=|(g_m(y_m)-g_m(y_0))(q_m-p_m)|> $$
\begin{equation}\label{eq3A}
> |(g_m(y_m)-g_m(y_0))q_m|=|(g_m(y_m)-g_m(y_0))q_m
+g_m(y_0)-g_m(y_0)|=
\end{equation}
$$=|g_m(y_0)-z_1^m|=\varepsilon_2^m\,.$$
Observe that
\begin{equation}\label{eq5B}
|g_m(y_0)- g_m(y_m)|+\varepsilon_1^m=|g_m(y_0)-
g_m(y_m)|+|g_m(y_m)-z_1^m|= |z_1^m-g_m(y_0)|=\varepsilon_2^m\,,
\end{equation}
consequently, $\varepsilon_1^m<\varepsilon_2^m.$ Now, we obtain
from~(\ref{eq3A}) that
\begin{equation}\label{eq4A}
|\gamma(1)-z_1^m|>\varepsilon_1^m\,.
\end{equation}

If $\gamma(0)\not\in S(z_1^m, \varepsilon_1^m),$ then we obtain by
(\ref{eq2A}) and (\ref{eq4A}) that $|\gamma|\cap B(z_1^m,
\varepsilon_1^m)\ne\varnothing\ne (D\setminus B(z_1^m,
\varepsilon_1^m))\cap|\gamma|.$ Thus, by \cite[Theorem~1.I, Ch.~5,
\S\, 46]{Ku} there exists $t_1\in (0, 1)$ such that $\gamma(t_1)\in
S(z_1^m, \varepsilon_1^m).$ Without loss of generality, we can
consider that $\gamma(t)\not\in B(z_1^m, \varepsilon_1^m)$ for $t\in
(t_1, 1).$ Put $\gamma_1:=\gamma|_{[t_1, 1]}.$

\medskip
From other hand, since $\varepsilon_1^m<\varepsilon_2^m$ and
$\gamma_1(t_1)\in S(z_1^m, \varepsilon_1^m)$, we obtain that
$|\gamma_1|\cap B(z_1^m, \varepsilon_2^m).$ By~(\ref{eq3A}), we
obtain that $(D\setminus B(z_1^m, \varepsilon_2^m))\cap
|\gamma_1|\ne\varnothing,$ so, by~\cite[Theorem~1.I, Ch.~5, \S\,
46]{Ku} there exists $t_2\in [t_1, 1)$ such that $\gamma_1(t_2)\in
S(z_1^m, \varepsilon_2^m).$ Without loss of generality, we can
consider that $\gamma_1(t)\in B(z_1^m, \varepsilon_2^m)$ for $t\in
(t_1, t_2).$ Put $\gamma_2:=\gamma|_{[t_1, t_2]}.$ Now
$\gamma>\gamma_2$ and $\gamma_2\in \Gamma(S(z_1^m, \varepsilon_1^m),
S(z_1^m, \varepsilon_2^m), A_m).$ Thus, (\ref{eq1A}) has been
proved.

\medskip
Put
$$\eta(t)= \left\{
\begin{array}{rr}
\frac{1}{\varepsilon_0}, & t\in [\varepsilon_1^m, \varepsilon_2^m],\\
0,  &  t\not\in [\varepsilon_1^m, \varepsilon_2^m]\,.
\end{array}
\right. $$
Observe that $\eta$ satisfies~(\ref{eq8B}) for
$r_1=\varepsilon_1^m,$ $r_2=\varepsilon_2^m.$ Indeed,
by~(\ref{eq13***}) and (\ref{eq5B}) we obtain that
$$r_1-r_2=\varepsilon_2^m-\varepsilon_1^m=|g_m(y_0)-z_1^m|-|g_m(y_m)-
z_1^m|=$$$$=|g_m(y_m)-g_m(y_0)|\geqslant\varepsilon_0\,.$$
Now,
$\int\limits_{\varepsilon_1^m}^{\varepsilon_2^m}\eta(t)dt=(1/\varepsilon_0)\cdot
(\varepsilon_2^m-\varepsilon_1^m)\geqslant 1.$

By the definition of ring $Q$-homeomorphism at the point $z_1^m$
and~(\ref{eq1A}), we obtain that
$$M(f_m(\Gamma_m))\leqslant M(f_m(\Gamma(S(z_1^m, \varepsilon_1^m),
S(z_1^m, \varepsilon_2^m), A_m)))\leqslant$$
\begin{equation}\label{eq14***}
\leqslant \frac{1}{\varepsilon_0^n}\int\limits_{D}
Q(x)\,dm(x):=c<\infty\,,
\end{equation}
as $Q\in L^1(D).$
On the other hand, $h(f_m(P_m))\geqslant h(y_m, f_m(z^m_1))
\geqslant (1/2)\cdot h(y_0, p_1)>0$ and $h(f_m(Q_m))\geqslant h(y_0,
f_m(z^m_2)) \geqslant (1/2)\cdot h(y_0, p_2)>0$ for large $m\in
{\Bbb N},$ where $h(f_m(Q_m))$ is defined in (\ref{eq9C}) for
$E:=f_m(Q_m).$ Moreover,
$$h(f_m(P_m), f_m(Q_m)):=\inf\limits_{x\in f_m(P_m), y\in f_m(Q_m)}h(x, y)\leqslant h(y_m, y_0)\rightarrow 0,\quad m\rightarrow
\infty\,.$$ By Lemma~\ref{lem2}
$$M(f_m(\Gamma_m))=M(f_m(P_m), f_m(Q_m), D^{\,\prime})\rightarrow\infty\,,\quad m\rightarrow\infty\,,$$
which contradicts the relation~(\ref{eq14***}). The contradiction
obtained above disproves the assumption in~(\ref{eq13***}). Theorem
has been proved.~$\Box$

\section{On behavior of mappings in the closure of domain}

Let us to turn to questions concerning the global behavior of
mappings. The following assertion indicates that, for sufficiently
good domains and mappings with condition~(\ref{eq2*!}), the image of
fixed continuum under mappings can not be close to the boundary of
the mapped domain, whenever Euclidean diameter of this continuum is
bounded from below (cf.~\cite[Theorems~21.13 and 21.14]{Va}).

\medskip
\begin{lemma}\label{lem3}
{\sl\, Let $D$ be a bounded domain in ${\Bbb R}^n,$ $n\geqslant 2,$
and let $D^{\,\prime}$ be a domain in $\overline{{\Bbb R}^n}.$
Suppose that $D$ is locally connected on $\overline{D},$
$D^{\,\prime}$ has weakly flat boundary, $Q\in L^1(D)$ and, besides
that, any connected component of $\partial D^{\,\prime}$ does not
degenerate to a point. Let $f_m:D\rightarrow D^{\,\prime}$ be a
sequence of ring $Q$-homeomorphisms in $D$ of $D$ onto
$D^{\,\prime}.$ Let $A\subset D$ be a continuum, and let $\delta>0$
be a number such that $h(f_m(A))\geqslant \delta>0$ for all
$m=1,2,\ldots ,$ where $h(f_m(A))$ is defined by~(\ref{eq9C}). Then
there exists $\delta_1>0$ such that
$$h(f_m(A),
\partial D^{\,\prime})>\delta_1>0\quad \forall\,\, m\in {\Bbb
N}\,,$$
where $h(f_m(A),
\partial D^{\,\prime})=\inf\limits_{x\in f_m(A), y\in \partial D^{\,\prime}}h(x, y).$}
\end{lemma}

\medskip
\begin{proof}
Since $D$ is bounded, and $f_m(D)=D^{\,\prime},$ $m=1,2,\ldots ,$ we
obtain that $\partial D^{\,\prime}\ne\varnothing.$ Thus, $h(f_m(A),
\partial D^{\,\prime})$ is well-defined.

\medskip
Assume the contrary. Now, for each $k\in {\Bbb N}$ there exists
$m=m_k:$ $h(f_{m_k}(A),
\partial D^{\,\prime})<1/k.$
Without loss of generality, we may assume that $m_k$ is increasing
sequence of numbers. Since $\overline{{\Bbb R}^n}$ is compact,
$\partial D^{\,\prime}$ is compact, as well. Note that $f_{m_k}(A)$
is a compact set as a continuous image of a compact set $A\subset D$
under the mapping $f_{m_k}.$ Now, there exist $x_k\in f_{m_k}(A)$
and $y_k\in
\partial D^{\,\prime}$ such that $h(f_{m_k}(A),
\partial D^{\,\prime})=h(x_k, y_k)<1/k$ (see~Figure~\ref{fig3}).
\begin{figure}[h]
\centerline{\includegraphics[scale=0.6]{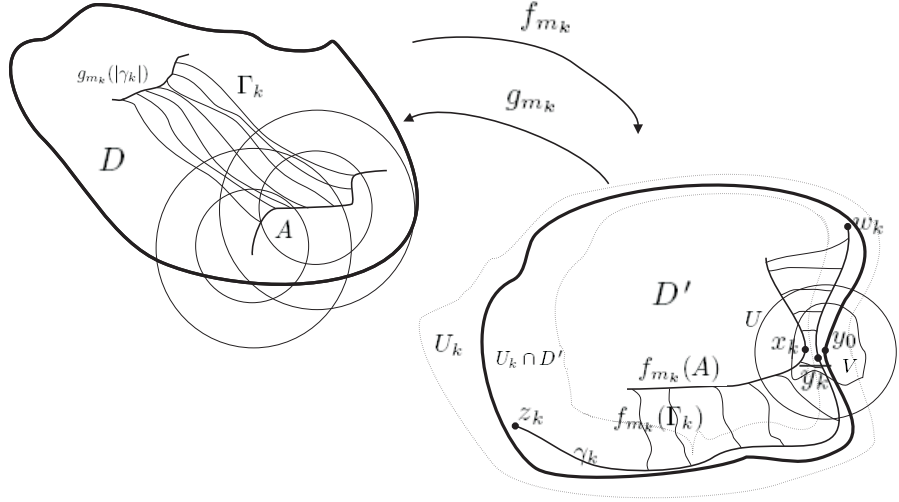}} \caption{To
the proof of Lemma~\ref{lem3}}\label{fig3}
\end{figure}
Since $\partial D^{\,\prime}$ is compact, we may assume that
$y_k\rightarrow y_0\in \partial D^{\,\prime},$ $k\rightarrow
\infty;$ then also
%
%\begin{equation}\label{eq7}
$$x_k\rightarrow y_0\in \partial D^{\,\prime},\quad k\rightarrow
\infty\,.$$
%\end{equation}
%
Let $K_0$ be a connected component of $\partial D^{\,\prime},$
containing $y_0.$ Obviously, $K_0$ is a continuum in
$\overline{{\Bbb R}^n}.$ Since $D^{\,\prime}$ has a weakly flat
boundary, the mapping $g_{m_k}:=f_{m_k}^{\,-1}$ extends to a
continuous mapping
$\overline{g}_{m_k}:\overline{D^{\,\prime}}\rightarrow \overline{D}$
for all $k\in {\Bbb N}$  (see~\cite[Theorem~3]{Sm}). Furthermore,
$\overline{g}_{m_k}$ is uniformly continuous on
$\overline{D^{\,\prime}},$ because $\overline{g}_{m_k}$ is a
continuous mapping on the compact set $\overline{D^{\,\prime}}.$ In
this case, for every $\varepsilon>0$ there is
$\delta_k=\delta_k(\varepsilon)<1/k$ such that
\begin{equation}\label{eq3}
|\overline{g}_{m_k}(x)-\overline{g}_{m_k}(x_0)|<\varepsilon \quad
\forall\,\, x,x_0\in \overline{D^{\,\prime}},\quad h(x,
x_0)<\delta_k\,, \quad \delta_k<1/k\,.
\end{equation}
Let $\varepsilon>0$ be such that
\begin{equation}\label{eq5}
\varepsilon<(1/2)\cdot {\rm  dist}\,(\partial D, A)\,.
\end{equation}
Denote $B_h(x_0, r)=\{x\in \overline{{\Bbb R}^n}: h(x, x_0)<r\}.$
Given $k\in {\Bbb N},$ we put
$$B_k:=\bigcup\limits_{x_0\in K_0}B_h(x_0, \delta_k)\,,\quad k\in {\Bbb N}\,.$$
Since $B_k$ is a neighborhood of $K_0,$ by~\cite[Lemma~2.2]{HK}
there exists a neighborhood $U_k$ of $K_0,$ such that $U_k\subset
B_k$ and $U_k\cap D^{\,\prime}$ is connected. Without loss of
generality, we may assume that $U_k$ is an open set, then $U_k\cap
D^{\,\prime}$ is also linearly connected
(see~\cite[Proposition~13.1]{MRSY}). Let $h(K_0)=m_0,$ where
$h(K_0)$ is defined in (\ref{eq9C}) for $E:=K_0.$ Now, we can find
$z_0, w_0\in K_0$ such that $h(K_0)=h(z_0, w_0)=m_0.$ Hence, we can
choose sequences $\overline{y_k}\in U_k\cap D^{\,\prime},$ $z_k\in
U_k\cap D^{\,\prime}$ and $w_k\in U_k\cap D^{\,\prime}$ such that
$z_k\rightarrow z_0,$ $\overline{y_k}\rightarrow y_0$ and
$w_k\rightarrow w_0$ as $k\rightarrow\infty.$ We may assume that
\begin{equation}\label{eq2}
h(z_k, w_k)>m_0/2,\quad \forall\,\, k\in {\Bbb N}\,.
\end{equation}
Since $U_k\cap D^{\,\prime}$ is path-connected, we can join the
points $z_k,$ $\overline{y_k}$ and $w_k$ sequentially by some path
$\gamma_k\in U_k\cap D^{\,\prime}.$ Let $|\gamma_k|$ be a locus of
$\gamma_k$ in $D^{\,\prime}.$ Now, $g_{m_k}(|\gamma_k|)$ is a
compact set in $D.$ If $x\in|\gamma_k|,$ then there is $x_0\in K_0:$
$x\in B(x_0, \delta_k).$ Put $\omega\in A\subset D.$ Since
$x\in|\gamma_k|$ and $x$ is an interior point of the domain
$D^{\,\prime},$ we write $g_{m_k}(x)$ instead of
$\overline{g}_{m_k}(x)$ in this case. By (\ref{eq3}) and (\ref{eq5})
and by triangle inequality, we obtain:
$$|g_{m_k}(x)-\omega|\geqslant
|\omega-\overline{g}_{m_k}(x_0)|-|\overline{g}_{m_k}(x_0)-g_{m_k}(x)|\geqslant$$
\begin{equation}\label{eq4}
\geqslant {\rm  dist}\,(\partial D, A)-(1/2)\cdot{\rm
dist}\,(\partial D, A)=(1/2)\cdot{\rm  dist}\,(\partial D,
A)>\varepsilon
\end{equation}
for sufficiently large $k\in {\Bbb N},$ ${\rm  dist}\,(\partial D,
A):=\inf\limits_{x\in \partial D, y\in A}|x-y|.$ Letting to $\inf$
in~(\ref{eq4}) over all $x\in |\gamma_k|$ and all $\omega\in A,$ we
obtain, that
\begin{equation}\label{eq6}
{\rm dist}\,(g_{m_k}(|\gamma_k|), A):=\inf\limits_{x\in
g_{m_k}(|\gamma_k|), y\in A}|x-y|>\varepsilon, \quad\forall\,\,
k=1,2,\ldots \,.
\end{equation}
Let us to cover the continuum $A$ by balls $B(x, \varepsilon/4),$
$x\in A.$ Since $A$ is the continuum, we can consider that $A\subset
\bigcup\limits_{i=1}^{M_0}B(x_i, \varepsilon/4),$ $x_i\in A,$
$i=1,2,\ldots, M_0,$ $1\leqslant M_0<\infty.$ By definition, $M_0$
depends only on $A,$ in particular, $M_0$ does non depend on $k.$
Putting
\begin{equation}\label{eq5C}
\Gamma_k:=\Gamma(A, g_{m_k}(|\gamma_k|), D)\,,
\end{equation}
we observe that
\begin{equation}\label{eq6C}
\Gamma_k=\bigcup\limits_{i=1}^{M_0}\Gamma_{ki}\,,
\end{equation}
where $\Gamma_{ki}$ consists of all paths $\gamma:[0, 1]\rightarrow
D$ from $\Gamma_k$ such that $\gamma(0)\in B(x_i, \varepsilon/4)$
and $\gamma(1)\in g_{m_k}(|\gamma_k|).$ Let us to show that
\begin{equation}\label{eq7C}
\Gamma_{ki}>\Gamma(S(x_i, \varepsilon/4), S(x_i, \varepsilon/2),
A(x_i, \varepsilon/4, \varepsilon/2))\,.
\end{equation}
Indeed, let $\gamma\in \Gamma_{ki},$ i.e., $\gamma:[0, 1]\rightarrow
D,$ $\gamma(0)\in B(x_i, \varepsilon/4),$ and $\gamma(1)\in
g_{m_k}(|\gamma_k|).$ By~(\ref{eq6}), $|\gamma|\cap B(x_i,
\varepsilon/4)\ne\varnothing\ne |\gamma|\cap (D\setminus B(x_i,
\varepsilon/4)).$ Thus, by \cite[Theorem~1.I, Ch.~5, \S\, 46]{Ku},
there exists $0<t_1<1$ with $\gamma(t_1)\in S(x_i, \varepsilon/4).$
We can consider that $\gamma(t)\not\in B(x_i, \varepsilon/4)$ for
$t>t_1.$ Put $\gamma_1:=\gamma|_{[t_1, 1]}.$ By~(\ref{eq6}),
$|\gamma_1|\cap B(x_i, \varepsilon/2)\ne\varnothing\ne
|\gamma_1|\cap (D\setminus B(x_i, \varepsilon/2)).$ Thus, by
\cite[Theorem~1.I, Ch.~5, \S\, 46]{Ku}, there exists $t_1<t_2<1$
with $\gamma(t_2)\in S(x_i, \varepsilon/2).$ We can consider that
$\gamma(t)\in B(x_i, \varepsilon/2)$ for $t<t_2.$ Put
$\gamma_2:=\gamma|_{[t_1, t_2]}.$ So, $\gamma_2$ is a subpath of
$\gamma,$ which belongs to $\Gamma(S(x_i, \varepsilon/4), S(x_i,
\varepsilon/2), A(x_i, \varepsilon/4, \varepsilon/2)).$ So, we have
proved~(\ref{eq7C}).
\medskip
Put
$$\eta(t)= \left\{
\begin{array}{rr}
4/\varepsilon, & t\in [\varepsilon/4, \varepsilon/2],\\
0,  &  t\not\in [\varepsilon/4, \varepsilon/2]\,.
\end{array}
\right. $$
Observe that $\eta$ satisfies~(\ref{eq8B}) for $r_1=\varepsilon/4,$
$r_2=\varepsilon/2.$ Now, by the definition of ring
$Q$-homeomorphism at $x_i$
\begin{equation}\label{eq8C}
M(f_{m_k}(\Gamma(S(x_i, \varepsilon/4), S(x_i, \varepsilon/2)),
A(x_i, \varepsilon/4, \varepsilon/2)))\leqslant
(4/\varepsilon)^n\cdot\Vert Q\Vert_1<c<\infty\,,
\end{equation}
where $c$ is some positive constant, and $\Vert Q\Vert_1$ is
$L_1$-norm of the function $Q$ in $D.$ By~(\ref{eq6C}), (\ref{eq7C})
and (\ref{eq8C}), using subadditivity of modulus, we obtain that
\begin{equation}\label{eq4B}
M(f_{m_k}(\Gamma_k))\leqslant
\frac{4^nM_0}{\varepsilon^n}\int\limits_DQ(x)\,dm(x)=c=c(\varepsilon,
Q)<\infty\,.
\end{equation}
\medskip
Let us to show that we obtain the contradiction of (\ref{eq4B}) with
weakly flatness of the boundary. Let $U:=B_h(y_0, r_0),$ where
$r_0>0,$ $r_0<\min\{\delta/4, m_0/4\},$ $\delta$ is a number from
the condition of the lemma, and $h(K_0)=m_0.$ (Here $h(K_0)$ denotes
the chordal diameter of a set $E=K_0,$ see~(\ref{eq9C})). Notice,
that $|\gamma_k|\cap U\ne\varnothing\ne |\gamma_k|\cap
(D^{\,\prime}\setminus U)$ for sufficiently large $k\in{\Bbb N},$
because the $h(|\gamma_k|)\geqslant m_0/2>m_0/4,$ $\overline{y_k}\in
|\gamma_k|$ and $\overline{y_k}\rightarrow y_0$ as
$k\rightarrow\infty.$ Similarly, $f_{m_k}(A)\cap U\ne\varnothing\ne
f_{m_k}(A)\cap (D^{\,\prime}\setminus U).$ Since $|\gamma_k|$ and
$f_{m_k}(A)$ are continua, we obtain that
\begin{equation}\label{eq8}
f_{m_k}(A)\cap \partial U\ne\varnothing, \quad|\gamma_k|\cap
\partial U\ne\varnothing\,,
\end{equation}
see~\cite[Theorem~1.I, Ch.~5, \S\, 46]{Ku}. Given $P>0,$ let
$V\subset U$ be a neighborhood of $y_0$ from the definition of a
weakly flat boundary. Now
\begin{equation}\label{eq9}
M(\Gamma(E, F, D^{\,\prime}))>P
\end{equation}
for any continua $E, F\subset D^{\,\prime}$ with $E\cap
\partial U\ne\varnothing\ne E\cap \partial V$ and $F\cap \partial
U\ne\varnothing\ne F\cap \partial V.$
Observe that
\begin{equation}\label{eq10}
f_{m_k}(A)\cap \partial V\ne\varnothing, \quad|\gamma_k|\cap
\partial V\ne\varnothing
\end{equation}
for sufficiently large $k\in {\Bbb N}.$ Indeed, $\overline{y_k}\in
|\gamma_k|,$ $x_k\in f_{m_k}(A),$ where $x_k,
\overline{y_k}\rightarrow y_0\in V$ as $k\rightarrow\infty.$
Therefore, $|\gamma_k|\cap V\ne\varnothing\ne f_{m_k}(A)\cap V$ for
large $k\in {\Bbb N}.$ Besides that, $h(V)\leqslant h(U)\leqslant
2r_0<m_0/2.$ By~(\ref{eq2}), $h(|\gamma_k|)>m_0/2,$ thus
$|\gamma_k|\cap (D^{\,\prime}\setminus V)\ne\varnothing.$ Therefore,
by~~\cite[Theorem~1.I, Ch.~5, \S\, 46]{Ku}, $|\gamma_k|\cap\partial
V\ne\varnothing.$. Similarly, $h(V)\leqslant h(U)\leqslant
2r_0<\delta/2.$ Since $h(f_{m_k}(A))>\delta,$ we obtain that
$f_{m_k}(A)\cap (D^{\,\prime}\setminus V)\ne\varnothing.$
By~\cite[Theorem~1.I, Ch.~5, \S\, 46]{Ku}, we have that
$f_{m_k}(A)\cap
\partial V\ne\varnothing.$ Now, (\ref{eq10})
is proved.

\medskip
By~(\ref{eq8}), (\ref{eq9}) and (\ref{eq10}), we obtain that
\begin{equation}\label{eq11}
M(\Gamma(f_{m_k}(A), |\gamma_k|, D^{\,\prime}))>P\,.
\end{equation}
Notice, that $\Gamma(f_{m_k}(A), |\gamma_k|,
D^{\,\prime})=f_{m_k}(\Gamma(A, g_{m_k}(|\gamma_k|),
D))=f_{m_k}(\Gamma_k).$ Thus, (\ref{eq11}) can be rewritten as
$$M(\Gamma(f_{m_k}(A), g_{m_k}(|\gamma_k|), D))=M(f_{m_k}(\Gamma_k))>P\,.$$
The relation obtained above contradicts (\ref{eq4B}). Thus, the
assumption ${\rm dist\,}(f_{m_k}(A),
\partial D^{\,\prime})<1/k$ was not true. The lemma is proved.~$\Box$
\end{proof}

\medskip
{\it Proof of Theorem~\ref{th2}}. Let $g\in {\frak S}_{\delta, A, Q
}(D, D^{\,\prime}).$ Since $D^{\,\prime}$ has a weakly flat
boundary, $g$ extends to a continuous mapping
$\overline{g}:\overline{D^{\,\prime}}\rightarrow \overline{D}$
(see~\cite[Theorem~3]{Sm}, cf.~\cite[Theorem~4.6]{MRSY}).

Let us to verify the equality
$\overline{g}(\overline{D^{\,\prime}})=\overline{D}.$ Indeed, by
definition,
$\overline{g}(\overline{D^{\,\prime}})\subset\overline{D}.$ It
remains to show the converse inclusion $\overline{D}\subset
\overline{g}(\overline{D^{\,\prime}}).$ Let $x_0\in \overline{D}.$
Now, we show that $x_0\in \overline{g}(\overline{D^{\,\prime}}).$ If
$x_0\in \overline{D},$ then either $x_0\in D,$ or $x_0\in
\partial D.$ If $x_0\in D,$ then there is nothing to prove, since by hypothesis
$\overline{g}(D^{\,\prime})=D.$  Let $x_0\in \partial D.$ Now, there
exist $x_k\in D$ and $y_k\in D^{\,\prime}$ such that
$x_k=\overline{g}(y_k)$ and $x_k\rightarrow x_0$ as
$k\rightarrow\infty.$ Since $\overline{D^{\,\prime}}$ is compact, we
may assume that $y_k\rightarrow y_0\in \overline{D^{\,\prime}}$ as
$k\rightarrow\infty.$ Since $f=g^{\,-1}$ is a homeomorphism, $y_0\in
\partial D^{\,\prime}.$ Since $\overline{g}^{\,-1}$ is continuous in $\overline{D^{\,\prime}},$
$\overline{g}(y_k)\rightarrow \overline{g}(y_0).$ However, in this
case, $\overline{g}(y_0)=x_0,$ because $\overline{g}(y_k)=x_k$ and
$x_k\rightarrow x_0$ as $k\rightarrow\infty.$ Thus, $x_0\in
\overline{g}(\overline{D^{\,\prime}}).$ The inclusion
$\overline{D}\subset \overline{g}(\overline{D^{\,\prime}})$ is
proved. Thus, $\overline{D}=\overline{g}(\overline{D^{\,\prime}}),$
as required.

\begin{figure}[h]
\centerline{\includegraphics[scale=0.6]{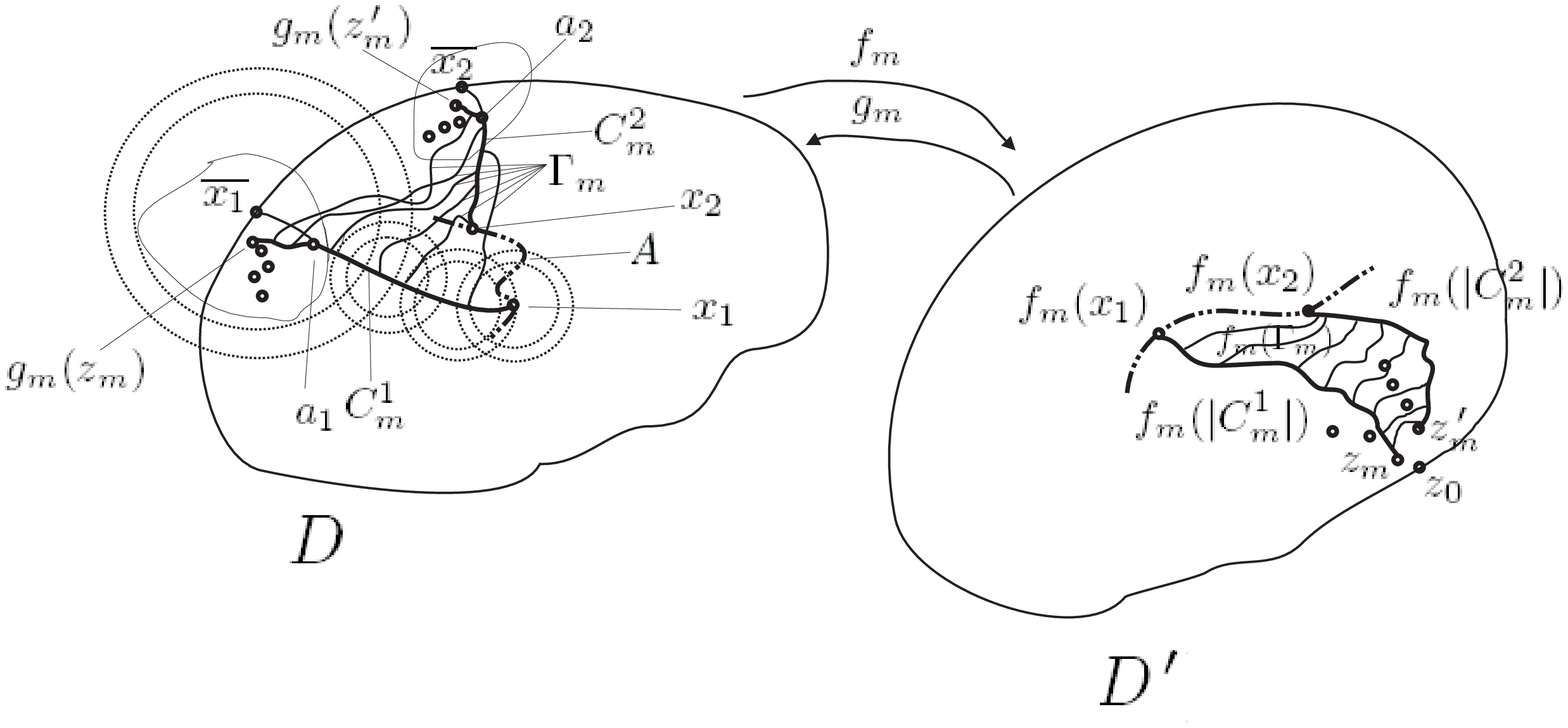}} \caption{To
the proof of Theorem~\ref{th2}}\label{fig4}
\end{figure}

The equicontinuity of ${\frak S}_{\delta, A, Q }(\overline{D},
\overline{D^{\,\prime}})$ in $D^{\,\prime}$ immediately follows from
Theorem~\ref{th1}. It remains to show that ${\frak S}_{\delta, A, Q
}(\overline{D}, \overline{D^{\,\prime}})$ is equicontinuous at
boundary points. We give the proof by contradiction. Now, we can
find a point $z_0\in \partial D^{\,\prime},$ a number
$\varepsilon_0>0$ and sequences $z_m\in \overline{D^{\,\prime}},$
$z_m\rightarrow z_0$ as $m\rightarrow\infty$ and $\overline{g}_m\in
{\frak S}_{\delta, A, Q }(\overline{D}, \overline{D^{\,\prime}})$
such that
\begin{equation}\label{eq12}
|\overline{g}_m(z_m)-\overline{g}_m(z_0)|\geqslant\varepsilon_0,\quad
m=1,2,\ldots .
\end{equation}
Put $g_m:=\overline{g}_m|_{D^{\,\prime}}.$ Since $g_m$ extends by
continuity to the boundary of $D^{\,\prime},$ we may assume that
$z_m\in D^{\,\prime}$ and, hence, $\overline{g}_m(z_m)=g_m(z_m).$ In
addition, there exists $z^{\,\prime}_m\in D^{\,\prime},$
$z^{\,\prime}_m\rightarrow z_0$ as $m\rightarrow\infty,$ such that
$|g_m(z^{\,\prime}_m)-\overline{g}_m(z_0)|\rightarrow 0$ as
$m\rightarrow\infty.$
%Тогда из (\ref{eq12}) вытекает, что
%%
%\begin{equation}\label{eq13}
%|g_m(z_m)-g_m(z^{\,\prime}_m)|\geqslant\varepsilon_0/2,\quad
%m\geqslant m_0\,.
%\end{equation}
%
Since $D$ is bounded, $\overline{D}$ is compact. Thus, we may assume
that $g_m(z_m)$ and $\overline{g}_m(z_0)$ are convergent sequences
as $m\rightarrow\infty.$ Assume that $g_m(z_m)\rightarrow
\overline{x_1}$ and $\overline{g}_m(z_0)\rightarrow \overline{x_2}$
as $m\rightarrow\infty.$ By continuity of the modulus in
(\ref{eq12}), $\overline{x_1}\ne \overline{x_2}.$ Besides that,
since homeomorphisms preserve a boundary, $\overline{x_2}\in\partial
D.$ Let $x_1$ and $x_2$ be arbitrary distinct points of the
continuum $A,$ none of which coincide with с $\overline{x_1}.$ By
Lemma~\ref{lem1} we can join points $x_1$ and $\overline{x_1}$ by
the path $\gamma_1:[0, 1]\rightarrow \overline{D},$ and points $x_2$
and $\overline{x_2}$ by the path $\gamma_2:[0, 1]\rightarrow
\overline{D}$ such that $|\gamma_1|\cap |\gamma_2|=\varnothing,$
$\gamma_i(t)\in D$ for all $t\in (0, 1),$ $i=1,2,$
$\gamma_1(0)=x_1,$ $\gamma_1(1)=\overline{x_1},$ $\gamma_2(0)=x_2$
and $\gamma_2(1)=\overline{x_2}.$ Since $D$ is locally connected on
$\partial D,$ there are neighborhoods $U_1$ and $U_2$ of
$\overline{x_1}$ and $\overline{x_2},$ whose closures do not
intersect, and $W_i:=D\cap U_i$ are path-connected sets. Without
loss of generality, we may assume that $\overline{U_1}\subset
B(\overline{x_1}, \delta_0)$ and
\begin{equation}\label{eq12C}
\overline{B(\overline{x_1},
\delta_0)}\cap|\gamma_2|=\varnothing=\overline{U_2}\cap|\gamma_1|\,,
\quad \overline{B(\overline{x_1}, \delta_0)}\cap
\overline{U_2}=\varnothing\,,
\end{equation}
$g_m(z_m)\in W_1$ and $g_m(z^{\,\prime}_m)\in W_2$ for each $m\in
{\Bbb N}.$ Let $a_1$ and $a_2$ be arbitrary points belonging to
$|\gamma_1|\cap W_1$ and $|\gamma_2|\cap W_2.$ Let $t_1, t_2$ be
such that $\gamma_1(t_1)=a_1$ and $\gamma_2(t_2)=a_2.$ We join $a_1$
and $g_m(z_m)$ by a path $\alpha_m:[t_1, 1]\rightarrow W_1$ such
that $\alpha_m(t_1)=a_1$ and $\alpha_m(1)=g_m(z_m).$ Similarly, we
join $a_2$ and $g_m(z^{\,\prime}_m)$ by a path $\beta_m:[t_2,
1]\rightarrow W_2,$ $\beta_m(t_2)=a_2$ and
$\beta_m(1)=g_m(z^{\,\prime}_m)$ (see Figure~\ref{fig4}).
\medskip
Set
$$C^1_m(t)=\quad\left\{
\begin{array}{rr}
\gamma_1(t), & t\in [0, t_1],\\
\alpha_m(t), & t\in [t_1, 1]\end{array} \right.\,,\qquad
C^2_m(t)=\quad\left\{
\begin{array}{rr}
\gamma_2(t), & t\in [0, t_2],\\
\beta_m(t), & t\in [t_2, 1]\end{array} \right.\,.$$
Denote, as usual, $|C^1_m|$ and $|C^2_m|$ are loci of paths $C^1_m$
and $C^2_m,$  respectively. Setting
$$l_0=\min\{{\rm dist}\,(|\gamma_1|,
|\gamma_2|), {\rm dist}\,(|\gamma_1|, U_2)\}\,,$$
we consider the covering $A_0:=\bigcup\limits_{x\in |\gamma_1|}B(x,
l_0/4).$ Since $|\gamma_1|$ is a compact, we can choose $1\leqslant
N_0<\infty$ and points $x_1,\ldots, x_{N_0}\in |\gamma_1|$ such that
$|\gamma_1|\subset B_0:=\bigcup\limits_{i=1}^{N_0}B(x_i, l_0/4).$
Now
$$|C^1_m|\subset U_1\cup |\gamma_1|\subset \overline{B(\overline{x_1}, \delta_0)}\cup \bigcup\limits_{i=1}^{N_0}B(x_i, l_0/4)\,.$$
Let $\Gamma_m$ be a family of paths connecting $|C^1_m|$ and
$|C^2_m|$ in $D.$ Now,
\begin{equation}\label{eq10C}
\Gamma_m=\bigcup\limits_{i=0}^{N_0}\Gamma_{mi}\,,
\end{equation}
where $\Gamma_{mi}$ consists of all paths $\gamma:[0, 1]\rightarrow
D$ with $\gamma(0)\in B(x_i, l_0/4)\cap |C^1_m|$ and $\gamma(1)\in
|C_2^m|$ for $1\leqslant i\leqslant N_0.$ Similarly, $\Gamma_{m0}$
consists of all paths $\gamma:[0, 1]\rightarrow D$ with
$\gamma(0)\in B(\overline{x_1}, \delta_0)\cap |C^1_m|$ and
$\gamma(1)\in |C_2^m|.$ By~(\ref{eq12C}) there exists
$\sigma_0>\delta_0>0$ such that
$$
\overline{B(\overline{x_1},
\sigma_0)}\cap|\gamma_2|=\varnothing=\overline{U_2}\cap|\gamma_1|\,,
\quad \overline{B(\overline{x_1}, \sigma_0)}\cap
\overline{U_2}=\varnothing\,.$$
Arguing similarly to proof of Lemma~\ref{lem3}, we can show that
$$\Gamma_{m0}>\Gamma(S(\overline{x_1}, \delta_0), S(\overline{x_1}, \sigma_0), A(\overline{x_1}, \delta_0,
\sigma_0))\,,$$
\begin{equation}\label{eq11C}
\Gamma_{mi}>\Gamma(S(x_i, l_0/4), S(x_i, l_0/2), A(x_i, l_0/4,
l_0/2))\,.
\end{equation}
Putting
$$\eta(t)= \left\{
\begin{array}{rr}
4/l_0, & t\in [l_0/4, l_0/2],\\
0,  &  t\not\in [l_0/4, l_0/2]\,,
\end{array}
\right. \qquad \eta_0(t)= \left\{
\begin{array}{rr}
1/(\sigma_0-\delta_0), & t\in [\delta_0, \sigma_0],\\
0,  &  t\not\in [\delta_0, \sigma_0]\,,
\end{array}
\right.$$
and $f_m:=g_m^{\,-1},$ we obtain by~(\ref{eq2*!}) that
$$
M(f_m(\Gamma(S(\overline{x_1}, \delta_0), S(\overline{x_1},
\sigma_0), A(\overline{x_1}, \delta_0, \sigma_0))))\leqslant
(1/(\sigma_0-\delta_0))^n\cdot\Vert Q\Vert_1<c_1<\infty\,,
$$
\begin{equation}\label{eq12D}
M(f_m(\Gamma(S(x_i, l_0/4), S(x_i, l_0/2), A(x_i, l_0/4,
l_0/2))))\leqslant (4/(l_0))^n\cdot\Vert Q\Vert_1<c_2<\infty\,,
\end{equation}
where $c_1$ and $c_1$ are some positive constants, not depending on
$m.$ We conclude from~(\ref{eq10C}), (\ref{eq11C}), (\ref{eq12D})
and subadditivity of modulus that
\begin{equation}\label{eq14A}
M(f_m(\Gamma_m))\leqslant (N_0/l_0^n+(1/(\sigma_0-\delta_0))^n)\Vert
Q\Vert_1:=c<\infty\,.
\end{equation}
From other hand, by Lemma~\ref{lem3}, there is a number $\delta_1>0$
such that $h(f_{m}(A), \partial D^{\,\prime})>\delta_1>0,$
$m=1,2,\ldots \,.$ Thus,
$$h(f_m(|C^1_m|))\geqslant h(z_m, f_m(x_1)) \geqslant
(1/2)\cdot h(f_m(A), \partial D^{\,\prime})>\delta_1/2\,,$$
\begin{equation}\label{eq14}
h(f_m(|C^2_m|))\geqslant h(z^{\,\prime}_m,f_m(x_2)) \geqslant
(1/2)\cdot h(f_m(A), \partial D^{\,\prime})>\delta_1/2
\end{equation}
for some $M_0\in {\Bbb N}$ and for all $m\geqslant M_0.$
%Кроме того, ${\rm dist}\,(f_m(|C^1_m|, f_m(|C^2_m|)\leqslant
%|z_m-z^{\,\prime}_m|\rightarrow 0$ по выбору $z_m$ и
%$z^{\,\prime}_m.$
%
Set $U:=B_h(z_0, r_0),$ where $0<r_0<\delta_1/4$ and $\delta_1$ is
from (\ref{eq14}). Notice, that $f_m(|C^1_m|)\cap U\ne\varnothing\ne
f_m(|C^1_m|)\cap (D^{\,\prime}\setminus U)$ for sufficiently large
$m\in{\Bbb N},$ because $h(f_m(|C^1_m|))\geqslant \delta_1/2$ and
$z_m\in f_m(|C^1_m|),$ $z_m\rightarrow z_0$ as $m\rightarrow\infty.$
Similarly, $f_m(|C^2_m|)\cap U\ne\varnothing\ne f_m(|C^2_m|)\cap
(D^{\,\prime}\setminus U).$ Since $f_m(|C^1_m|)$ and $f_m(|C^2_m|)$
are continua,
\begin{equation}\label{eq8A}
f_m(|C^1_m|)\cap \partial U\ne\varnothing, \quad f_m(|C^2_m|)\cap
\partial U\ne\varnothing\,,
\end{equation}
see, e.g.,~\cite[Theorem~1.I, Ch.~5, \S\, 46]{Ku}. Since $\partial
D^{\,\prime}$ is weakly flat, given $P>0,$ there exists a
neighborhood $V\subset U$ of $z_0$ such that
\begin{equation}\label{eq9A}
M(\Gamma(E, F, D^{\,\prime}))>P
\end{equation}
for any continua $E, F\subset D^{\,\prime}$ with $E\cap
\partial U\ne\varnothing\ne E\cap \partial V$ и $F\cap \partial
U\ne\varnothing\ne F\cap \partial V.$ Observe that
\begin{equation}\label{eq10A}
f_m(|C^1_m|)\cap \partial V\ne\varnothing, \quad f_m(|C^2_m|)\cap
\partial V\ne\varnothing\end{equation}
for sufficiently large $m\in {\Bbb N}.$

Indeed, let $z_m\in f_m(|C^1_m|),$ $z^{\,\prime}_m\in f_m(|C^2_m|),$
where $z_m, z^{\,\prime}_m\rightarrow z_0\in V$ as
$m\rightarrow\infty.$ Now, $f_m(|C^1_m|)\cap V\ne\varnothing\ne
f_m(|C^2_m|)\cap V$ for sufficiently large $m\in {\Bbb N}.$ In
addition, $h(V)\leqslant h(U)\leqslant 2r_0<\delta_1/2.$ Besides
that, by (\ref{eq14}) we obtain that $h(f_m(|C^1_m|))>\delta_1/2.$
Thus, $f_m(|C^1_m|)\cap (D^{\,\prime}\setminus V)\ne\varnothing$
and, consequently, $f_m(|C^1_m|)\cap\partial V\ne\varnothing$
(see~\cite[Theorem~1.I, Ch.~5, \S\,~46]{Ku}). Similarly,
$h(V)\leqslant h(U)\leqslant 2r_0<\delta_1/2.$ By~(\ref{eq14})
$h(f_m(|C^2_m|))>\delta,$ thus $f_m(|C^2_m|)\cap
(D^{\,\prime}\setminus V)\ne\varnothing.$ By~\cite[Theorem~1.I,
Ch.~5, \S\, 46]{Ku} we have, that $f_m(|C^1_m|)\cap\partial
V\ne\varnothing.$ Thus, (\ref{eq10A}) is proved.

\medskip
By~(\ref{eq8A}), (\ref{eq9A}) and (\ref{eq10A}), we obtain that
%
%\begin{equation}\label{eq11}
$$M(f_m(\Gamma_m))=M(\Gamma(f_m(|C^1_m|), f_m(|C^2_m|),
D^{\,\prime}))>P\,,$$
%\end{equation}
%
which contradicts~(\ref{eq14A}). The contradiction obtained above
disproves the assumption made in~(\ref{eq12}). The theorem is
proved.~$\Box$

\section{Some examples}

We begin with a simple example of mappings on the complex plane.

\medskip
{\bf Example~1.} As known, the linear-fractional automorphisms of
the unit disk ${\Bbb D}\subset{\Bbb C}$ onto itself can be written
by the formula $f(z)=e^{i\theta}\frac{z-a}{1-\overline{a}z},$ $z\in
{\Bbb D},$ $a\in{\Bbb C},$ $|a|<1,$ $\theta\in [0, 2\pi).$ These
mappings $f$ are ring 1-homeomorphisms; all the conditions of
Theorem~\ref{th2} are satisfied, except the condition
$h(f(A))\geqslant\delta,$ which, in general, does not hold

If, for instance, $\theta=0,$ $a=1/n,$ $n=1,2,\ldots,$ then
$f_n(z)=\frac{z-1/n}{1-z/n}=\frac{nz-1}{n-z}.$ Let $A=[0, 1/2].$ Now
$f_n(0)=-1/n\rightarrow 0$ and $f_n(1/2)=\frac{n-2}{2n-1}\rightarrow
1/2$ as $n\rightarrow\infty.$ Thus, $f_n$ satisfies the condition
$h(f_n(A))\geqslant\delta$ for $\delta=1/4.$ We obtain, that
$f_n^{\,-1}(z)=\frac{z+1/n}{1+z/n}$ and, hence, $f_n^{\,-1}$
converge uniformly to $f^{\,-1}(z)\equiv z.$ Thus, the sequence
$f_n^{\,-1}(z)$ is equicontinuous in $\overline{{\Bbb D}}.$

\medskip
Now, put
$f^{\,-1}_n(z)=\frac{z-(n-1)/n}{1-z(n-1)/n}=\frac{nz-n+1}{n-nz+1}.$
It is easy to see, that $f^{\,-1}_n$ converges locally uniformly to
$-1$ inside of ${\Bbb D},$ whenever $f^{\,-1}_n(1)=1.$ Now, we
conclude that $f^{\,-1}_n$ is not equicontinuous at 1. In this case,
$f_n(z)=\frac{z+(n-1)/n}{1+z(n-1)/n}$ and the condition
$h(f_n(A))\geqslant\delta$ does not hold for any $\delta>0$
by~Theorem~\ref{th2}.

\medskip
Thus, under the hypotheses of Theorem~\ref{th2}, we can not refuse
from the additional requirement $h(f(A))\geqslant\delta,$ in
general.

\medskip
{\bf Example~2.} Let $p\geqslant 1$ be a number, such that
$n/p(n-1)<1.$ Put $\alpha\in (0, n/p(n-1)).$ We define a sequence of
mappings $f_m:{\Bbb B}^n\rightarrow B(0, 2)$ of ${\Bbb B}^n$ onto
the $B(0, 2)$ in the following way:
$$f_m(x)\,=\,\left
\{\begin{array}{rr} \frac{1+|x|^{\alpha}}{|x|}\cdot x\,, & 1/m\leqslant|x|\leqslant 1, \\
\frac{1+(1/m)^{\alpha}}{(1/m)}\cdot x\,, & 0<|x|< 1/m \ .
\end{array}\right.
$$
Notice, that $f_m$ satisfies (\ref{eq2*!}) for
$Q=\left(\frac{1+|x|^{\,\alpha}}{\alpha
|x|^{\,\alpha}}\right)^{n-1}\in L^1({\Bbb B}^n)$ at every $x_0\in
\overline{{\Bbb B}^n},$ see~\cite[proof of Theorem~7.1]{Sev$_3$}.
By~\cite[Lemma~4.3]{Vu}, $B(0, 2)$ has a weakly flat boundary.
Observe that $f_m$ fixes an infinite number of points of the unit
ball for all $m\geqslant 2.$

\medskip
By Theorem~\ref{th2}, the family $\frak{G}=\{g_m\}_{m=1}^{\infty},$
$g_m:=f_m^{\,-1},$ is equicontinuous in $\overline{B(0, 2)}.$

\medskip
Observe that the ''inverse'' family $\frak F=\{f_m\}_{m=1}^{\infty}$
is not equicontinuous in ${\Bbb B}^n.$ Indeed,
$|f_m(x_m)-f(0)|=1+1/m\not\rightarrow 0$ as $m\rightarrow\infty,$
where $|x_m|=1/m$).

\medskip
{\it The family $\frak G$ contains an infinite number of mappings
$g_{m_k}:=f^{\,-1}_{m_k},$ $f_{m_k}\in \frak F,$ that do not satisfy
the relation (\ref{eq2*!}) with $Q\in L^1$}. Indeed, otherwise,
by~Theorem~\ref{th1} ''the inverse''  to $\frak G$ family $\frak F$
must be equicontinuous in ${\Bbb B}^n.$

%=================Список литературы====================
%\end{fulltext}

\medskip
\medskip
{\bf \noindent Evgeny Sevost'yanov, Sergei Skvortsov} \\
Zhytomyr Ivan Franko State University,  \\
40 Bol'shaya Berdichevskaya Str., 10 008  Zhytomyr, UKRAINE \\
Phone: +38 -- (066) -- 959 50 34, \\
Email: esevostyanov2009@gmail.com

\end{document}